\documentclass[reqno,10pt]{amsart}
\usepackage{hyperref} 

\usepackage{geometry}
\usepackage[utf8]{inputenc}
\usepackage{enumerate}
\usepackage{amstext,amsmath,amsthm,amsfonts,amssymb,amscd}
\usepackage{latexsym,mathrsfs,dsfont,euscript}
\usepackage{mathbbol}
\usepackage{mathtools}
\usepackage{esint}
\usepackage{array,float}
\usepackage{xcolor,colortbl}
\usepackage{fullpage}
\usepackage{soul}
\usepackage{subcaption}
\usepackage{enumitem}
\usepackage{wasysym}
\usepackage{ragged2e}
\usepackage{tabularx}

\hypersetup{
  colorlinks,%
  citecolor=blue,%
    filecolor=red,%
    linkcolor=red,%
    urlcolor=blue
}

\usepackage{graphicx}
\newtheorem{theorem}{Theorem}[section]
\newtheorem{proposition}[theorem]{Proposition}
\newtheorem{lemma}[theorem]{Lemma}
\newtheorem{corollary}[theorem]{Corollary}
\theoremstyle{definition}
\newtheorem{definition}[theorem]{Definition}

\theoremstyle{remark}

\numberwithin{equation}{section}

\allowdisplaybreaks

\begin{document}

%%%%%%%%%%%%% Metadata %%%%%%%%%%%%%%%%%%%%%%%%%%%%%%%%%%%
%% Title
\title[titulo corto]{A measure theoretic approach to Lipschitz regularity and its Haar type wavelet analysis}

%%    Information for authors

\author[]{Hugo Aimar}
\email{haimar@santafe-conicet.gov.ar}
%%%
\author[]{Juliana Boasso}
\email{jboasso@santafe-conicet.gov.ar}
%%%
\date{\today}
%\email{xxx@santafe-conicet.gov.ar}
%%

%%  Palabras claves y codigos

%ejemplos
%53B12  	Differential geometric aspects of statistical manifolds and information geometry


\begin{abstract}
    In this note, we extend the characterization of dyadic Lipschitz regularity of functions to non-atomic probability spaces, using generalized Haar systems.
\end{abstract}

\maketitle

\section{Introduction}\label{sec:Section1}
    As it is well known from the results obtained by M. Holschneider and P. Tchamitchian, \cite{HT90} and \cite{HT91}, see also \cite{D92}, $\alpha$-Lipschitz regularity of functions defined on Euclidean spaces is completely characterized in terms of the behavior of the continuous wavelet transform. These results hold true for very general wavelet families. The basic feature of this characterization is given by the power law on the scale, uniformly in the position parameter, for the size of the projection of the function under analysis. The classical Haar function $h$ in $\mathbb{R}$ shows that such a characterization is not possible in terms of the behavior of the discrete wavelet coefficients with respect to the standard orthonormal basis generated by $h.$ Nevertheless, the results in \cite{AAG24} show that some hidden regularity is still reflected by a power law in the scale, uniformly in the position parameter, for the Haar wavelet coefficients. In fact, this behavior characterizes Lipschitz regularity with respect to the distance induced by the measure of the dyadic intervals. On the other hand, this metric is only given by the dyadic family itself and the measure space supporting the dyadic sets. The usual dyadic families in $\mathbb{R}^n$ and even the generalizations built on spaces of homogeneous type, see \cite{Chr90} and \cite{Dav88}, are strongly related to the underlying metric structure. In particular, they are differentiation bases in the sense defined in \cite{deG81}. This differentiation property leads to the construction of orthonormal Haar wavelet-type bases for the space $L^2.$ Non-metric guided dyadic families can be constructed in very general settings. Also, specially simple constructions in the plane $\mathbb{R}^2$ can be designed in order to produce quantitative parameters measuring particular textures in images and time series. The wavelet-type functions provided by these general dyadic families will usually not produce orthonormal bases for the space of square-integrable functions. Even so, as we aim to show in this note, there is no need for the basis character of the Haar system induced by a dyadic family in order to obtain a characterization through power laws of the Haar coefficients of the Lipschitz type regularity of a function with respect to the metric induced by the dyadic family. The exponents of these power laws would serve as indicators of textures and diverse features of signals and images. The Hurst exponent (see \cite{H51}) used in hydrology as a measure of persistence for time series can be computed using wavelets. Let us observe that when there is a metric in the space guiding the construction of the dyadic families, the results in \cite{ABN12} show that the dyadic families are suitable to recover the metric Lipschitz character of a function. It is worthy to mention at this point the construction of regular wavelet bases given in \cite{AH13} (see also \cite{AB23}) in spaces of homogeneous type. 
    
    In this note we work on probability non-atomic spaces, without metric, but equipped with a family of dyadic sets that does not need to be a differentiation basis for the Lebesgue spaces of the setting.

    \vspace{1cm}
    %\medskip
    The paper is organized as follows. Section~\ref{sec:Section2} is devoted to introducing the dyadic families on a non-atomic probability space and some of their basic properties. In Section~\ref{sec:Section3}, we review the definition of pseudo-metrics produced by dyadic systems in measure spaces. In Section~\ref{sec:Section4}, we introduce the notions of Lipschitz regularity with respect to a given dyadic family that we shall consider later. Section~\ref{sec:Section5} is devoted to reviewing the construction of the Haar wavelet families provided by general dyadic systems and to proving some of their basic properties. The main result, i.e., the characterization of dyadic Lipschitz regularity in terms of the generalized Haar coefficients, is proved in Section~\ref{sec:Section6}. The main results are contained in Theorems~\ref{teo61}, \ref{teo63}, and \ref{teo66}. Finally, Section~\ref{sec:Section7} is devoted to illustrating, in the case of elementary 2-D images, the ability of different dyadic families to reflect particular textures.

\section{Dyadic systems in probability spaces} \label{sec:Section2}
    \hspace{0.6cm}Let $(\mathbb{X},\mathcal{F},\mu)$ be a finite positive measure space. For simplicity and without loosing generality we shall assume that $\mu(\mathbb{X})=1.$ Or, in other words, that $(\mathbb{X},\mathcal{F},\mu)$ is a probability space. Let us recall that the measurable sets, the elements of $\mathcal{F},$ admit a notion of distance given by the measure of the symmetric difference $\mu(E\vartriangle F).$ In this setting the notions of density and separability are given by the general context of metric spaces. As it can be seen, for example, in \cite{Hal74}, every finite measure space $(\mathbb{X},\mathcal{F},\mu)$ that is separable and non-atomic, is naturally isomorphic to the interval $[0,1]$ with the Lebesgue measure. This fact allows the construction on any non-atomic separable probability space $(\mathbb{X},\mathcal{F},\mu)$ of a diversity of dyadic systems, as introduced in the next definition. It is worthy noticing that this isomorphism is only relevant to answer the questions regarding the existence of dyadic families satisfying the required properties in a general setting. Usually, for example, in the case of the analysis of images, the dyadic families are crafted by partitioning in several ways the frame of the given image. We shall see some examples in Section~\ref{sec:Section7}.
    
    \begin{definition}\label{def21}
        We shall say that a family $\mathcal{D}$ in a probability space $(\mathbb{X},\mathcal{F},\mu)$ is a \textbf{dyadic family} with \textbf{inheritance coefficient} $B\geq 2$ in $(\mathbb{X},\mathcal{F},\mu)$ if

        \begin{enumerate}[label=\textbf{\thetheorem.\alph*}, leftmargin=*, align=left]
            \item $\mathcal{D}\subset\mathcal{F};$ \label{def21a}
            \item $\mathcal{D}=\bigcup_{j \geq 0} \mathcal{D}^j,$ with $\mathcal{D}^0=\lbrace\mathbb{X}\rbrace;$ and for $j\geq 1,$ $\mathcal{D}^j=\big\lbrace Q_{j,k}: k\in\mathcal{K}_j\rbrace, $ $\#\lbrace\mathcal{K}_j\rbrace=K_j,\:$ $\bigcup_{k\in\mathcal{K}_j} Q_{j,k}=\mathbb{X}$ and $ Q_{j,k}\cap Q_{j,l}=\emptyset$ for $k\neq l$ both in $\mathcal{K}_j;$\label{def21b}
            \item $\mu_{j,k}\doteq\mu(Q_{j,k})>0,$ for all $j\geq 0$ and $k\in\mathcal{K}_j;$\label{def21c}
            \item for each $j\geq 1$ and each $k\in\mathcal{K}_j$ there is a unique  $l\in\mathcal{K}_{j-1}$ such that $Q_{j,k}\subsetneq Q_{j-1,l},$ we briefly say that $Q_{j-1,l}$ is the first ancestor of $Q_{j,k};$\label{def21d}
            \item if $Q$ is the first ancestor of $\tilde{Q},$ both in $\mathcal{D},$ then $\mu(Q)\leq B\mu(\tilde{Q}).$\label{def21e}
        \end{enumerate}
    \end{definition}
    The next statement provides some elementary but relevant properties of a dyadic family.
    \begin{lemma}\label{lemma22}
        Let $\mathcal{D}$ be a dyadic family with inheritance coefficient $B.$ Then
        \begin{enumerate}[label=\textbf{\thetheorem.\arabic*}, leftmargin=*, align=left]
            \item  $\bigcup_{\tilde{Q}\in\mathcal{O}(Q)}\tilde{Q} = Q,$ for every $Q\in\mathcal{D}$ where $\mathcal{O}$ denotes the offspring of $Q,$ in other words $\mathcal{O}(Q)=\lbrace\tilde{Q}: Q \text{ is the first ancestor of }\tilde{Q}\rbrace;$ \label{lem221}
            \item given any $Q\in\mathcal{D},$ the offspring of $Q$ has at least two elements and at most $B;$\label{lem222}
            \item given $j$ and $l$ positive integers, $k\in\mathcal{K}_j$ and $m\in\mathcal{K}_{j+l}$ with $Q_{j+l,m}\subset Q_{j,k},$ then $\mu_{j+l,m}\leq \gamma^l \mu_{j,k},$ with $0<\gamma=1-\frac{1}{B}<1.$\label{lem223}
        \end{enumerate}
        \begin{proof}
            To prove \ref{lem221} notice that from \ref{def21d} we have $\bigcup_{\tilde{Q}\in\mathcal{O}(Q)}\tilde{Q} \subset Q.$ Assume $Q=Q_{j,k}.$ Since $\mathcal{D}^{j+1}$ is a disjoint cover of $\mathbb{X},$ then $Q\subset\bigcup_{l\in\mathcal{K}_{j+1}}Q_{j+1,l}.$ If $l\in\mathcal{K}_{j+1}$ is such that $Q_{j+1,l}\cap Q \neq \emptyset,$ then again from \ref{def21d} we must have that $Q_{j+1,l}\subset Q.$ Or in other words $Q_{j+1,l}\in \mathcal{O}(Q).$ So that $Q\subset \bigcup_{\tilde{Q}\in\mathcal{O}(Q)}\tilde{Q}.$ In order to prove \ref{lem222}, let us first observe that, from \ref{def21d}, $\#\lbrace\mathcal{O}(Q)\rbrace \geq 2.$ In order to prove the upper bound for $\#\lbrace\mathcal{O}(Q)\rbrace$ write the measure of $Q$ as the sum of the measures of its offspring and use \ref{def21e}
            $$\mu(Q)=\mu\bigg(\bigcup_{\tilde{Q}\in\mathcal{O}(Q)}\tilde{Q}\bigg)=\sum_{\tilde{Q}\in\mathcal{O}(Q)}\mu(\tilde{Q}) \geq \sum_{\tilde{Q}\in\mathcal{O}(Q)}\frac{1}{B}\mu(Q) = \frac{1}{B}\#\lbrace\mathcal{O}(Q)\rbrace \mu(Q),$$ which implies that $\#\lbrace\mathcal{O}(Q)\rbrace\leq B,$ as stated. Let us finally prove \ref{lem223}. Suppose first that $l=1,$ since $Q_{j,k}$ is the first ancestor of $Q_{j+1,l}$ and from \ref{lem222}, $Q_{j+1,m}$ has at least one brother, say  $Q_{j+1,p},$ then from \ref{def21e} we get $\mu_{j,k}\geq \mu_{j+1,m}+\mu_{j+1,p}\geq \mu_{j+1,m}+\frac{1}{B}\mu_{j,k}.$ In other words $\mu_{j+1,m}\leq \big(1-\frac{1}{B}\big)\mu_{j,k}=\gamma \mu_{j,k}.$ So, the result follows inductively in $l.$
            
        \end{proof}
    \end{lemma}

    Of course, the basic examples of dyadic systems in probability spaces are the classical, usually metric driven constructions such as those in the interval $[0,1),$ where $\mathcal{D}^j=\lbrace I_{j,k}=[k 2^{-j},(k+1)2^{-j}): k=0,\dots, 2^j-1\rbrace$, with $j$ a non-negative integer.
    Or in a cube $Q=[0,1)^m$ of $\mathbb{R}^m,$ with $\mathcal{D}^j=\lbrace Q_{j,\mathbf{k}}=\prod_{i=1}^m I_{j,k_i}: \mathbf{k}=(k_1,\dots, k_m), k_i=0,\dots, 2^j-1; i=1,\dots,m \rbrace. $
    Or even particular rectangles in $Q$ like $\mathcal{D}^j=\lbrace R_{j,k}=\big[k_1 4^{-j}, (k_1+1)4^{-j}\big)\times\prod_{i=2}^m I_{j,k_i}: k_1=0,\dots,4^j-1; k_i=0,\dots, 2^j-1; i=2,\dots,m \rbrace.$ All of them are balls with respect to some metric in $\mathbb{X}.$ The first are balls with respect to distances which are equivalent to the Euclidean. The last are balls associated to parabolic distances in $\mathbb{R}^m.$ Being so, these families share several extra properties aside from \ref{def21a} to \ref{def21e} such as that of being differentiation bases in the Lebesgue sense. Nevertheless, our approach, purely measure theoretic, lacks of metric structure and allows a large diversity of partitions. For example, for $Q=[0,1)^2=\mathbb{X},$ an acceptable dyadic family for the square is given by $\mathcal{D}^j=\lbrace R_{j,k}=I_{j,k}\times [0,1): k=0,\dots,2^j-1 \rbrace, j\geq 0,$ which satisfies properties \ref{def21a} to \ref{def21e} with $B=2,$ but it is not a differentiation basis in the sense of \cite{deG81}.

    The basic fact is that these non-differentiating dyadic families still provide Haar type wavelets, which in general will only be orthonormal systems in $L^2,$ but not orthonormal bases. These systems are still enough to characterize a certain regularity of Lipschitz type for real functions defined in $\mathbb{X}.$ The Lipschitz character is defined through a pseudo metric generated by the same family $\mathcal{D}$ that are on the basis of the construction of the wavelets. Of course, when the ``cubes'' in $\mathcal{D}$ are driven by a metric or quasi-metric, like in \cite{Chr90}, we recover more classical results of characterization of regularity through the behavior of the Haar coefficients; see \cite{ABN12}.

\section{Pseudo metrization of probability spaces through dyadic families} \label{sec:Section3}
    \hspace{0.6cm}The aim of this brief section is to introduce the pseudo metric provided by a dyadic family on a probability space and to prove some of its elementary properties.
    \begin{proposition}\label{propo31}
        Let $\mathcal{D}$ be a dyadic family with inheritance coefficient $B$ in the probability space $(\mathbb{X},\mathcal{F},\mu).$ Then the function $$\delta_\mathcal{D}:\mathbb{X}\times\mathbb{X}\to \mathbb{R}^+\cup\lbrace 0 \rbrace$$ given by $$\delta_\mathcal{D}(x,y)=\inf\lbrace\mu(Q):x,y \in Q;Q\in\mathcal{D}\rbrace$$ satisfies the following properties.
        \begin{enumerate}[label=\textbf{\thetheorem.\arabic*}, leftmargin=*, align=left]
            \item  $\delta_\mathcal{D}(x,x)=0$ for every $x\in\mathbb{X};$ \label{propo311}
            \item  $\delta_\mathcal{D}(x,y)=\delta_\mathcal{D}(y,x),$ for every $x$ and $y$ in $\mathbb{X};$ \label{propo312}
            \item  $\delta_\mathcal{D}(x,z) \leq \sup\lbrace\delta_\mathcal{D}(x,y),\delta_\mathcal{D}(y,z)\rbrace$ for every choice of $x,$ $y$ and $z$ in $\mathbb{X}.$\label{propo313}
        \end{enumerate}
        
        \begin{proof}
            Given $x\in\mathbb{X}$ and $j$ a non-negative integer, from \ref{def21b} we have one and only one $k=k(j,x)\in \mathcal{K}_j$ such that $x\in Q_{j,k(j,x)}.$ Now, from \ref{lem223} we have $\mu(Q_{j,k(j,x)})\leq\gamma^j\mu(Q_{0,1})=\gamma^j.$ Therefore, $\mu(Q_{j,k(j,x)})$ approaches 0 as $j$ grows to infinity. Let $x\in\mathbb{X},$ then there are elements $Q\in\mathcal{D}$ of measure as small as we wish containing $x,$ so that $\delta_{\mathcal{D}}(x,x)=0.$ The second property, \ref{propo312}, follows directly from the symmetry of the condition that defines $\delta_{\mathcal{D}}.$ Finally, to prove \ref{propo313}, let $Q_1$ and $Q_2$ be the two dyadic cubes in $\mathcal{D}$ that realize the distances $\delta_\mathcal{D}(x,y)$ and $\delta_\mathcal{D}(y,z)$ respectively. Then, $\delta_\mathcal{D}(x,y)=\mu(Q_1)$ with $x,y\in Q_1$ and $\delta_\mathcal{D}(y,z)=\mu(Q_2)$ with $y,z\in Q_2.$ 
            Since $y\in Q_1\cap Q_2,$ it follows that $Q_1\subset Q_2$ or $Q_2\subset Q_1.$ Suppose  $Q_1\subset Q_2,$ then $x$ and $z$ are both in $Q_2.$ This implies that $$\delta_\mathcal{D}(x,z) \leq \mu(Q_2)=\sup \lbrace \mu(Q_1),\mu(Q_2)\rbrace=\sup \lbrace \delta_\mathcal{D}(x,y), \delta_\mathcal{D}(y,z)\rbrace.$$
            
        \end{proof}
    \end{proposition}

    It is worthy mention at this point that the reliability condition ($\delta_\mathcal{D}(x,y)=0$ implies $x=y$) for a metric does not hold in general in our setting. This is the case of the dyadic family, with $B=2,$ in the set $\mathbb{X}=[0,1)^2$ given by  $\mathcal{D}^j=\lbrace R_{j,k}=I_{j,k}\times [0,1): k=0,\dots, 2^j-1 \rbrace,$ $j\geq 0$ and, as before,  $I_{j,k}=[k 2^{-j}, (k+1)2^{-j}).$ With the usual terminology, we say that the $\delta_\mathcal{D}$ introduced in Proposition~\ref{propo31} is a pseudo metric or, more precisely, a pseudo-ultra metric. The next statement contains a simple but basic fact; the $\delta_\mathcal{D}$-balls are precisely the sets in $\mathcal{D}.$ As usual we shall consider in $\mathbb{X}$ the open $\delta_\mathcal{D}$ ball centered at $x\in\mathbb{X}$ with radius $r>0$ as the set $\mathcal{B}_{\delta_{\mathcal{D}}}(x,r)=\lbrace y\in\mathbb{X}:\delta_{\mathcal{D}}(x,y)<r\rbrace.$

    \begin{proposition}\label{propo32}
        For every $x\in\mathbb{X}$ and every $r>0,$ $$\mathcal{B}_{\delta_{\mathcal{D}}}(x,r)=Q(x,r),$$ where $Q(x,r)$ is the largest dyadic cube in $\mathcal{D}$ containing $x$ with measure strictly less than $r.$ So that, $$\mathcal{D}=\lbrace \mathcal{B}_{\delta_{\mathcal{D}}}(x,r): x\in\mathbb{X}, r>0\rbrace.$$
        \begin{proof}
            Notice that from the definition of $\delta_{\mathcal{D}}$ we have that $y\in\mathcal{B}_{\delta_{\mathcal{D}}}(x,r)$ if and only if there exists $Q\in\mathcal{D}$ containing both $x$ and $y$ with $\mu(Q)<r.$ Since given two cubes $Q$ and $Q'$ belonging to $\mathcal{D},$ both containing $x$ and $y,$ we have $Q\subset Q'$ or $Q'\subset Q,$ the result is clear.
            
        \end{proof}
    \end{proposition}
    
    Notice that for $r>1$ the $\delta_{\mathcal{D}}$-balls are always the whole space $\mathbb{X}.$ The following Proposition shows that the space $(\mathbb{X},\delta_{\mathcal{D}},\mu)$ is Ahlfors regular with exponent equal to one.
    \begin{proposition}\label{propo33}
        For $0<r<1,$ we have that $$\frac{r}{B}\leq \mu\big(\mathcal{B}_{\delta_{\mathcal{D}}}(x,r)\big)<r,$$ for every $x\in\mathbb{X}.$
        \begin{proof}
            From \ref{lem223}, we have that given $x,$ with the notation in the proof of Proposition~\ref{propo31}, the sequence $\lbrace\mu_{j,k(j,x)}:j\geq 0\rbrace$ is strictly decreasing to zero and $\mu_{0,1}=1.$ Then, given $0<r<1,$ there exists one and only one $j\geq 0$ such that $\mu_{j+1,k(j+1,x)}<r\leq \mu_{j,k(j,x)}.$ Hence $\mathcal{B}_{\delta_{\mathcal{D}}}(x,r)=Q(x,r)=Q_{j+1,k(j+1,x)}.$ So that, from \ref{def21e},
            \begin{align*}
                \frac{r}{B}\leq\frac{1}{B}\mu_{j,k(j,x)}=\frac{1}{B}\mu\big(Q_{j,k(j,x)}\big)&\leq \mu\big(Q_{j+1,k(j+1,x)}\big)\\ &=\mu\big(\mathcal{B}_{\delta_{\mathcal{D}}}(x,r)\big)= \mu_{j+1,k(j+1,x)}<r.
            \end{align*}
        
        \end{proof}
    \end{proposition}

    In a general metric space $(\mathbb{X},d),$ the distance from a subset $A$ of $\mathbb{X}$ to another subset $B$ of $\mathbb{X}$ is defined by $d(A,B)=\inf\lbrace d(x,y): x\in A,y\in B\rbrace.$ In our current dyadic setting, with $d=\delta_{\mathcal{D}},$ when the subsets $A$ and $B$ are themselves dyadic cubes, we have an explicit way to compute their distances.

    \begin{proposition}\label{propo34}
        Let $Q_1$ and $Q_2$ be two disjoint dyadic cubes in $\mathcal{D}.$ Then $$\delta_{\mathcal{D}}(Q_1,Q_2) = \inf\lbrace\mu(Q): Q\in\mathcal{D}\text{ and }Q_1\cup Q_2 \subset Q\rbrace.$$
        \begin{proof}
            Since $\delta_{\mathcal{D}}(Q_1,Q_2)=\inf\lbrace\delta_{\mathcal{D}}(x,y): x\in Q_1\text{ and } y \in Q_2\rbrace,$ we have to consider $\delta_{\mathcal{D}}(x,y)=\inf\lbrace\mu(Q): Q\in\mathcal{D};x\in Q\text{ and } y \in Q\rbrace$ for every choice of $x\in Q_1$ and $y \in Q_2.$ So that $\delta(Q_1,Q_2)=\inf_{\substack{x\in Q_1\\ y\in Q_2}}\inf_{\substack{x,y\in Q\\ Q\in \mathcal{D}}}\mu(Q).$ Since $Q_1\cap Q_2 = \emptyset,$ the smallest $Q\in\mathcal{D}$ containing some $x\in Q_1$ and some $y\in Q_2$ is the same as the smallest cube in $\mathcal{D}$ containing $Q_1$ and $Q_2.$ In fact, since $Q$ contains $x,$ then $Q\cap Q_1\neq \emptyset,$ if it would happen that $Q_1\not\subset Q,$ then it should happen that $Q\subset Q_1.$ But in this case, since $y\in Q\subset Q_1$ we have that $y\in Q_1\cap Q_2,$ which contradicts the hypothesis $Q_1\cap Q_2 = \emptyset.$ Then $\delta(Q_1,Q_2) =\inf_{\substack{Q_1\subset Q\\ Q_2\subset Q \\ Q\in \mathcal{D}}}\mu(Q).$ 
            
        \end{proof}
    \end{proposition}
    
\section{Lipschitz and generalized Lipschitz functions in $(\mathbb{X},\delta_{\mathcal{D}})$} \label{sec:Section4}
    \hspace{0.6cm}Once a metric or pseudo metric is given on a set $\mathbb{X},$ we have function spaces generated by regularity conditions provided by the metric. The best known are those given by power laws of the distance. The Lipschitz or Hölder functions.
    \begin{definition}\label{def41}
        Let $\mathcal{D}$ be a dyadic family with inheritance coefficient $B$ in the probability space $(\mathbb{X},\mathcal{F},\mu).$ Let $\alpha$ be a given positive real number. We say that a real valued function $f$ defined on $\mathbb{X}$ is of \textbf{class $\alpha$-Lipschitz with respect to the dyadic family $\mathcal{D},$ }if there exists some positive constant $A$ such that the inequality $$|f(x)-f(y)|\leq A\delta_\mathcal{D}^\alpha(x,y)$$ holds for every $x$ and $y$ in $\mathbb{X}.$ The infimum of the constants $A$ is called the semi-norm of $f$ and the vector space of $\alpha$-Lipschitz functions with respect to $\mathcal{D}$ is denoted by $Lip_{\delta_\mathcal{D}}(\alpha).$
    \end{definition}
     Let us observe at this point that, in Euclidean instances of our general setting, for example when $\mathbb{X}=[0,1)$ and $\mathcal{D}$ is the standard dyadic family, we have that $|x-y|\leq \delta_\mathcal{D}(x,y).$ Hence every standard (Euclidean) Lipschitz function is also a Lipschitz function, with the same exponent, with respect to $\mathcal{D}.$ In the general case, we have that the indicator function of $Q\in\mathcal{D}$ belongs to $Lip_{\delta_\mathcal{D}}(\alpha)$ for every $\alpha$ with semi-norm given by $\big(\mu(Q)\big)^{-\alpha}.$

     Two reasons leads us to consider a weaker form of the $\alpha$-Lipschitz condition with respect to a general dyadic family. The first is theoretical. Our characterization by Haar type wavelet systems without the differentiation property in $\mathcal{D}$ will be based in mean values. The second has practical purposes for applications where digital time series or images are given as mean values of some unknown continuous signal. The following result provides the key for the next definition.
     \begin{proposition}\label{propo42}
         Let $f$ be of class $\alpha$-Lipschitz with respect to the dyadic family $\mathcal{D}$ in $(\mathbb{X},\mathcal{D},\mu).$ Given $Q\in\mathcal{D},$ let $f_Q$ to denote the mean value of $f$ on $Q.$ Then the inequality $$\big|f_{Q_{j,m}}-f_{Q_{j,n}}\big|\leq A \delta_\mathcal{D}^\alpha(Q_{j,m},Q_{j,n})$$ holds for every $j\geq 0$ and every $m$ and $n$ on $\mathcal{K}_j.$
         \begin{proof}
             When $m=n$ there is nothing to prove. Assume that $n\neq m$ both in $\mathcal{K}_j$ and write
             \begin{align*}
                  \big|f_{Q_{j,m}}-f_{Q_{j,n}}\big| &\leq \frac{1}{\mu(Q_{j,m})}\int_{Q_{j,m}}\big| f(x) - f_{Q_{j,n}}\big| \:d\mu(x) \\ &= \frac{1}{\mu(Q_{j,m})}\int_{Q_{j,m}}\Big| f(x) -  \frac{1}{\mu(Q_{j,n})}\int_{Q_{j,n}}f(y)\:d\mu(y)\Big| \:d\mu(x) \\ &\leq \frac{1}{\mu(Q_{j,m})}\frac{1}{\mu(Q_{j,n})}\int_{Q_{j,m}}\int_{Q_{j,n}} |f(x)-f(y)|\:d\mu(y)d\mu(x) \\ &\leq \frac{A}{\mu(Q_{j,m})\mu(Q_{j,n})}\int_{Q_{j,m}}\int_{Q_{j,n}}\delta_\mathcal{D}^\alpha(x,y)\:d\mu(y)d\mu(x).
             \end{align*}
             Notice that, in the last integral $x\in Q_{j,m},$ $y\in Q_{j,n}$ and $Q_{j,m}\cap Q_{j,n}=\emptyset$ because $m\neq n.$ So that, as we proved in Proposition~\ref{propo34}, we have that $\delta_\mathcal{D}(x,y)$ is constantly equal to $\delta_\mathcal{D}(Q_{j,m},Q_{j,n})$ and the desired inequality is proved.
             
         \end{proof}
     \end{proposition}
     \begin{definition}\label{def43}
         Let $(\mathbb{X},\mathcal{D},\mu)$ and $\mathcal{D}$ be as before. A function $f$ defined and integrable on $\mathbb{X}$ is said to be of \textbf{pixelated $\alpha$-Lipschitz class with respect to $\mathcal{D}$} if, for some positive constant $A,$ the inequality $$\big|f_{Q_{j,m}}-f_{Q_{j,n}}\big|\leq A \delta_\mathcal{D}^\alpha(Q_{j,m},Q_{j,n})$$ holds for every $j\geq 0$ and every choice of $m$ and $n$ on $\mathcal{K}_j.$ We shall write $\Lambda_{\delta_\mathcal{D}}(\alpha)$ to denote the space of these functions and $|f|_{\Lambda_{\delta_\mathcal{D}}(\alpha)}$ to denote the infimum of those constants $A.$
     \end{definition}
    Proposition~\ref{propo42} shows that $\alpha$-Lipschitz implies pixelated $\alpha$-Lipschitz. In our general setting, when $\mathcal{D}$ is not a differentiation basis for $L^1(\mathbb{X},\mathcal{F},\mu)$ the converse is generally false. 

   In concrete applications to the analysis of images, there is only a finite number of scales available. This fact leads to the following definition.
    \begin{definition}\label{def44}
        Let $(\mathbb{X},\mathcal{D},\mu),$ $\mathcal{D}$ and $f$ be as in Definition~\ref{def43}. We shall say that $f$ is of \textbf{discrete $\alpha$-Lipschitz type up to level $J> 0$ with respect to $\mathcal{D}$ with bound $A$} if the inequality $$\big|f_{Q_{j,m}}-f_{Q_{j,n}}\big|\leq A \delta_\mathcal{D}^\alpha(Q_{j,m},Q_{j,n})$$ holds for every $j\in\lbrace0,1,\dots,J\rbrace$ and every $m$ and $n$ in $\mathcal{K}_j.$ We denote the space of these functions $f$ by $\Lambda_{\delta_\mathcal{D}}^J(\alpha)$ and we shall write $|f|_{\Lambda_{\delta_\mathcal{D}}^J(\alpha)}$ to denote the infimum of the constants $A.$
    \end{definition}
    
    \section{Haar type wavelets associated to dyadic families} \label{sec:Section5}
         \hspace{0.6 cm}Let $\mathcal{D}$ be a dyadic family with inheritance coefficient $B\geq 2$ in $(\mathbb{X},\mathcal{F},\mu).$ On the bases of the properties of $\mathcal{D},$ a multiresolution type analysis can be given in $L^2(\mathbb{X},\mathcal{F},\mu)$ by defining, for $j\geq 0$
         $$\mathcal{V}_j=\lbrace f: \mathbb{X}\to \mathbb{R} \text{ such that } f|_{Q_{j,k}}\text{ is constant for every } k\in\mathcal{K}_j\rbrace.$$ Notice that each $\mathcal{V}_j$ is a closed subspace of $L^2(\mathbb{X},\mathcal{F},\mu).$ Moreover $\mathcal{V}_j$ is a closed subspace of $\mathcal{V}_{j+1}$ with the $L^2(\mathbb{X},\mathcal{F},\mu)$ norm. Notice also that $\bigcap_{j\geq 0}\mathcal{V}_j=\mathcal{V}_0$ is the one dimensional space of all constant functions. Under our general assumptions, the density of $\bigcup_{j\geq 0}\mathcal{V}_j$ fails to be true.

         Given $Q\in \mathcal{D},$ the offspring $\mathcal{O}(Q)$ can be conveniently indexed by their scale and by the position index in that scale. Precisely, given $j\geq 0$ and $k\in \mathcal{K}_j,$ let $\mathcal{K}_{j+1,k}=\lbrace l \in \mathcal{K}_{j+1}:Q_{j+1,l}\subset Q_{j,k}\rbrace$ be the set of indices of the offspring of $Q_{j,k}.$ We shall also consider the auxiliary subspaces of $L^2(\mathbb{X},\mathcal{F},\mu)$ defined by $$\mathcal{V}_{j,k}=\lbrace f:\mathbb{X}\to \mathbb{R}/ f|_{Q_{j+1,l}} \text{ is constant }\forall \:l\in\mathcal{K}_{j+1,k}\text{ and }f(x)=0 \text{ for }x\notin Q_{j,k}\rbrace.$$
        Notice that from \ref{lem222}, $2\leq \dim \mathcal{V}_{j,k}\leq B,$ instead $\dim \mathcal{V}_{j}$ tends to infinity as $j$ grows.
    
        All the above construction allows to build a Haar type family of functions in $\mathbb{X}$ which will generally fail to be an orthonormal bases for $L^2(\mathbb{X},\mathcal{F},\mu),$ but it will suffice to characterize Lipschitz regularity as described in Section 4.
        \begin{proposition}\label{propo51}
            Let $\mathcal{D}$ be a dyadic family in $(\mathbb{X},\mathcal{F},\mu)$ with inheritance coefficient $B.$ Let $\mathcal{V}_{j,k}$ as above. Then 
            \begin{enumerate}[label=\textbf{\thetheorem.\arabic*}, leftmargin=*, align=left]
            \item for each $j\geq 0$ and each $k\in\mathcal{K}_j$ there exists an orthonormal basis of $\mathcal{V}_{j,k}$ given by $\lbrace (\mu_{j,k})^{-\frac{1}{2}}\chi_{Q_{j,k}}\rbrace \cup \mathcal{H}_{j,k},$ where $\mathcal{H}_{j,k}=\lbrace \psi_{j,k}^\lambda:\lambda=1,\dots, \dim \mathcal{V}_{j,k}-1\rbrace;$ \label{propo511}
            \item the system $\mathscr{H}_\mathcal{D}=\bigcup_{j\geq 0} \bigcup_{k\in\mathcal{K}_j}\mathcal{H}_{j,k}$ is orthonormal in $L^2(\mathbb{X},\mathcal{F},\mu).$ We shall say that $\mathscr{H}_\mathcal{D}$ is a Haar system generated by $\mathcal{D};$ \label{propo512}
            \item if $\mathcal{D}$ is a differentiation basis of $L^1(\mathbb{X},\mathcal{F},\mu),$ in the sense that for $f \in L^1(\mathbb{X},\mathcal{F},\mu),$ $f_{Q_{j,k(j,x)}}\to f(x)$ as $ j\to\infty$ for almost every $x,$ then the family $\lbrace 1\rbrace \cup \mathscr{H}_\mathcal{D}$ is an orthonormal basis for $L^2(\mathbb{X},\mathcal{F},\mu).$ \label{propo513} 
            \end{enumerate}
            \begin{proof}
                In order to prove \ref{propo511}, take $\mathcal{K}'_{j+1,k}$ to be any subset of $\mathcal{K}_{j+1,k}$ with $\#\lbrace \mathcal{K}'_{j+1,k}\rbrace= \dim \mathcal{V}_{j,k}-1 = \#\lbrace \mathcal{K}_{j+1,k}\rbrace-1.$ Then, the set of $\dim \mathcal{V}_{j,k}$ elements given by $\lbrace\chi_{Q_{j,k}}\rbrace \cup \lbrace \chi_{Q_{j+1,l}}:l\in\mathcal{K}'_{j+1,l}\rbrace$ is an algebraic basis of $\mathcal{V}_{j,k},$ because they are linearly independent. Since $\mathcal{V}_{j,k}$ inherits the Hilbert structure of $L^2(\mathbb{X},\mathcal{F},\mu),$ we can apply the Gramm-Schmidt algorithm starting with $\chi_{Q_{j,k}}.$ This procedure gives an orthonormal basis of $\mathcal{V}_{j,k}$ of $\dim\mathcal{V}_{j,k}$ elements, the first being $(\mu_{j,k})^{-\frac{1}{2}}\chi_{Q_{j,k}}(x).$ The remaining $\dim\mathcal{V}_{j,k}-1$ elements of this basis, which being orthogonal to $\chi_{Q_{j,k}}$ necessarily will have zero integral, are denoted by $\mathcal{H}_{j,k}=\lbrace \psi_{j,k}^\lambda:\lambda=1,\dots, \dim \mathcal{V}_{j,k}-1\rbrace.$ Of course, $\langle \psi_{j,k}^\lambda,\psi_{j,k}^{\tilde{\lambda}}\rangle=\delta_{\lambda,\tilde{\lambda}},$ with $\delta_{\lambda,\tilde{\lambda}}$ be the Kronecker delta. Let us now prove \ref{propo512}. From the construction of $\lbrace \psi_{j,k}^\lambda:j\geq 0; k\in\mathcal{K}_j;\lambda=1,\dots, \dim \mathcal{V}_{j,k}-1\rbrace$ we certainly have that $\|\psi_{j,k}^\lambda\|_2=1$ for every choice of $j,$ $k$ and $\lambda.$ So that we only need to show that $\langle \psi_{j,k}^\lambda,\psi_{j',k'}^{\lambda'}\rangle=\int_\mathbb{X} \psi_{j,k}^\lambda (x)\psi_{j',k'}^{\lambda'}(x)\: d\mu(x)=0$ when $j\neq j'$ or $k\neq k'$ or $\lambda\neq\lambda'.$ When $j=j',$ $k=k'$ and $\lambda\neq\lambda',$ the result follows from the construction of $\mathcal{H}_{j,k}.$ We only have to consider the case $j\neq j' $ and $Q_{j,k}\cap Q_{j',k'}\neq \emptyset.$ Assume that $j>j',$ then the support of $\psi_{j,k}^\lambda$ is contained in some cube of the offspring of $Q_{j',k'}.$ But since $\psi_{j',k'}^{\lambda'}$ is constant on each cube of the offspring of $Q_{j',k'}$ and the integral of $\psi_{j,k}^\lambda$ vanishes, we have that
                $$\int_\mathbb{X}\psi_{j',k'}^{\lambda'}(x)\psi_{j,k}^\lambda(x)\:d\mu(x)=C\int_\mathbb{X}\psi_{j,k}^\lambda(x)\:d\mu(x)=0.$$ Property \ref{propo513} follows, as usual, from the fact that under the hypothesis of differentiation of $\mathcal{D},$ we have that $\bigcup_{j\geq 0}\mathcal{V}_{j}$ is dense in $L^2(\mathbb{X},\mathcal{F},\mu).$
                
            \end{proof}
        \end{proposition}
    \section{The main result} \label{sec:Section6}
        \hspace{0.6 cm}In this section, we search for a characterization of dyadic $\alpha$-Lipschitz regularity of a function in $\mathbb{X}$ in terms of a power law on the scales for the Haar system defined by $\mathcal{D}.$ At this point it is worth noticing that the one-dimensional resemblance of our results is due to the fact that, as we proved in Proposition~\ref{propo33}, the space $(\mathbb{X},\mathcal{F},\mu)$ is 1-Ahlfors regular. So that when dealing with classical spaces such as $\mathbb{R}^d$ with the Lebesgue measure and the standard dyadic cubes, $\delta_{\mathcal{D}}(x,y)$ is the diameter of the smallest dyadic cube containing $x$ and $y,$ raised to the power $d.$ The first result is the simplest one and shows that the Lipschitz character of $f$ entails a power law of the Haar coefficients in terms of the measures of their dyadic supports. For this first result, we actually do not need the existence of an inheritance coefficient $B\geq 2$ and properties \ref{def21a} to \ref{def21d} of the dyadic family are sufficient. 
        \begin{theorem}\label{teo61}
            Let $\mathcal{D}$ be a dyadic family in $(\mathbb{X},\mathcal{F},\mu)$ that satisfies properties \ref{def21a} to \ref{def21d}. Let $\alpha>0$ and $f\in Lip_{\delta_\mathcal{D}}(\alpha)$ be given. Then, the inequality $$ \big|\langle f,\psi_{j,k}^\lambda\rangle\big|\leq |f|_{Lip_{\delta_\mathcal{D}}(\alpha)} (\mu_{j,k})^{\alpha+\frac{1}{2}}$$ holds for every $j\geq 0$, every $k\in\mathcal{K}_j$ and every $\lambda=1,2,\dots,\dim\mathcal{V}_{j,k}-1.$
            \begin{proof}
                Fix $j,k$ and $\lambda.$ Since $\psi_{j,k}^\lambda$ has a vanishing mean, we may write 
                \begin{align*}
                     \langle f,\psi_{j,k}^\lambda\rangle &= \int_\mathbb{X} f(x) \psi_{j,k}^\lambda(x) \: d\mu(x) \\ &= \int_\mathbb{X} \big(f(x)-f_{Q_{j,k}}\big)\psi_{j,k}^\lambda(x)\:d\mu(x) \\ &= \int_{Q_{j,k}} \frac{1}{\mu_{j,k}}\int_{Q_{j,k}} \big(f(x)-f(y)\big)\:\psi_{j,k}^\lambda(x)\:d\mu(y)\:d\mu(x).      
                \end{align*}
                Taking absolute values in the last identity and using the fact that $f\in Lip_{\delta_\mathcal{D}},$ we see that
                \begin{align*}
                    \big|\langle f,\psi_{j,k}^\lambda\rangle\big| &\leq \frac{1}{\mu_{j,k}}\int_{Q_{j,k}} \int_{Q_{j,k}} \big|f(x)-f(y)\big|\:d\mu(y)\:|\psi_{j,k}^\lambda(x)|\:d\mu(x) \\ &\leq |f|_{Lip_{\delta_\mathcal{D}}(\alpha)}\frac{1}{\mu_{j,k}}\int_{Q_{j,k}} \int_{Q_{j,k}} \delta_\mathcal{D}^\alpha(x,y) \:|\psi_{j,k}^\lambda(x)|\: d\mu(y)\:d\mu(x) \\ &\leq |f|_{Lip_{\delta_\mathcal{D}}(\alpha)}\frac{1}{\mu_{j,k}}\int_{x\in Q_{j,k}} \int_{y\in Q_{j,k}} \big(\mu(Q_{j,k})\big)^\alpha \:|\psi_{j,k}^\lambda(x)|\: d\mu(y)\:d\mu(x) \\ &= |f|_{Lip_{\delta_\mathcal{D}}(\alpha)}(\mu_{j,k})^\alpha \int_{x\in Q_{j,k}} |\psi_{j,k}^\lambda(x)|\:d\mu(x) \\ &\leq |f|_{Lip_{\delta_\mathcal{D}}(\alpha)}(\mu_{j,k})^\alpha\| \psi_{j,k}^\lambda \|_{2}\|\chi_{Q_{j,k}}\|_{2} \\ &= |f|_{Lip_{\delta_\mathcal{D}}(\alpha)}(\mu_{j,k})^{\alpha+\frac{1}{2}},
                \end{align*}
                 as desired.
                
            \end{proof}
        \end{theorem}
        In order to get a version of the above result for functions $f$ in the pixelated and discrete classes $\Lambda_{\delta_\mathcal{D}}(\alpha)$ and $\Lambda_{\delta_\mathcal{D}}^J(\alpha)$, defined in \ref{def43} and \ref{def44}, let us first prove the following result.
        \begin{lemma}\label{lemma62}
            Let $\mathcal{D}$ be a dyadic family in $(\mathbb{X},\mathcal{F},\mu)$ that satisfies properties \ref{def21a} to \ref{def21d}. Let $\alpha>0$ be given and $J$ a positive integer larger than two. Given $f\in L^1(\mathbb{X},\mathcal{F},\mu),$ define $g(x)=\sum_{m\in\mathcal{K}_J}f_{Q_{J,m}}\chi_{Q_{J,m}}(x).$ Then, the identity $$\int_{Q_{j,k}} f \psi_{j,k}^\lambda\:d\mu=\int_{Q_{j,k}} g \psi_{j,k}^\lambda\:d\mu$$ holds for every $j$ such that $0\leq j<J-1,$ every $k\in\mathcal{K}_j$ and every $\lambda=1,2,\dots,\dim\mathcal{V}_{j,k}-1.$
            \begin{proof}
                The result can be directly seen from the projections on the spaces $\mathcal{V}_j.$ Nevertheless we include the following direct computation based on the particular structure of the wavelets $\psi_{j,k}^\lambda.$ Since $\psi_{j,k}^\lambda\in\mathcal{V}_{j,k},$ then $\psi_{j,k}^\lambda$ can be written as $$\psi_{j,k}^\lambda(x)=\sum_{\lbrace l:Q_{j+1,l}\subset Q_{j,k}\rbrace}\beta_{j,k,\lambda,l}\:\chi_{Q_{j+1,l}}(x).$$ Then
                \begin{align*}
                    \int_{Q_{j,k}}f\psi_{j,k}^\lambda\:d\mu &= \int_{Q_{j,k}}f\Bigg(\sum_{\lbrace l:Q_{j+1,l}\subset Q_{j,k}\rbrace}\beta_{j,k,\lambda,l}\:\chi_{Q_{j+1,l}}\Bigg)\:d\mu \\ &= \sum_{\lbrace l:Q_{j+1,l}\subset Q_{j,k}\rbrace}\beta_{j,k,\lambda,l}\int_{Q_{j+1,l}}f\:d\mu \\ &= \sum_{\lbrace l:Q_{j+1,l}\subset Q_{j,k}\rbrace}\beta_{j,k,\lambda,l}\sum_{\lbrace m:Q_{J,m}\subset Q_{j+1,l}\rbrace}\int_{Q_{J,m}}f\:d\mu \\ &= \sum_{\lbrace l:Q_{j+1,l}\subset Q_{j,k}\rbrace}\beta_{j,k,\lambda,l}\sum_{\lbrace m:Q_{J,m}\subset Q_{j+1,l}\rbrace}f_{Q_{J,m}}\int_{Q_{j,k}}\chi_{Q_{J,m}}\:d\mu \\ &= \sum_{\lbrace l:Q_{j+1,l}\subset Q_{j,k}\rbrace}\beta_{j,k,\lambda,l}\int_{Q_{j+1,l}}g\:d\mu \\ &= \int_{Q_{j,k}}g\Bigg(\sum_{\lbrace l:Q_{j+1,l}\subset Q_{j,k}\rbrace}\beta_{j,k,\lambda,l}\:\chi_{Q_{j+1,l}}\Bigg)\:d\mu \\ &= \int_{Q_{j,k}}g\psi_{j,k}^\lambda\:d\mu.
                \end{align*}
                
            \end{proof}
        \end{lemma}
        \begin{theorem}\label{teo63}
            Let $\mathcal{D}$ be a dyadic family in $(\mathbb{X},\mathcal{F},\mu)$ that satisfies properties \ref{def21a} to \ref{def21d}. Let $\alpha>0$ be given. Let $f\in\Lambda_{\delta_\mathcal{D}}(\alpha)$ (resp. $f\in\Lambda_{\delta_\mathcal{D}}^J(\alpha)$). Then, the inequality $$ \big|\langle f,\psi_{j,k}^\lambda\rangle\big|\leq |f|_{\Lambda_{\delta_\mathcal{D}}(\alpha)}(\mu_{j,k})^{\alpha+\frac{1}{2}}, \: (\text{resp. }\big|\langle f,\psi_{j,k}^\lambda\rangle\big|\leq |f|_{\Lambda_{\delta_\mathcal{D}}^J(\alpha)}(\mu_{j,k})^{\alpha+\frac{1}{2}})$$ holds for every $j\geq 0$ (resp. $0\leq j\leq J-1$), every $k\in\mathcal{K}_j$ and every $\lambda=1,2,\dots,\dim\mathcal{V}_{j,k}-1.$
            \begin{proof}
                For large $J$ and $0\leq j< J-1,$ applying Lemma \ref{lemma62}, we have $$\langle f,\psi_{j,k}^\lambda\rangle=\langle g,\psi_{j,k}^\lambda\rangle = \int_{x\in Q_{j,k}}\frac{1}{\mu_{j,k}}\int_{y\in Q_{j,k}}\big(g(x)-g(y)\big)\psi_{j,k}^\lambda(x)\:d\mu(y)d\mu(x),$$ so that $$\big| \langle f,\psi_{j,k}^\lambda\rangle\big|\leq \frac{1}{\mu_{j,k}}\int_{x\in Q_{j,k}}\int_{y\in Q_{j,k}}|g(x)-g(y)| \big|\psi_{j,k}^\lambda(x)\big|\:d\mu(y)d\mu(x).$$ On the other hand, from the definition of $g$ we have that $|g(x)-g(y)|$ vanishes if $x$ and $y$ belong to $Q_{J,m}$ whereas $|g(x)-g(y)|=\big| f_{Q_{J,m}}-f_{Q_{J,n}}\big|$ if $x\in Q_{J,m}$ and $y\in Q_{J,n}.$ From the hypothesis on $f$ we have that $$|g(x)-g(y)|\leq |f|_{\Lambda_{\delta_\mathcal{D}}(\alpha)}\delta_{\mathcal{D}}^\alpha(Q_{J,m},Q_{J,n}).$$ Since both $Q_{J,m}$ and $Q_{J,n}$ are sub-cubes of $Q_{j,k},$ we then have that $$|g(x)-g(y)|\leq |f|_{\Lambda_{\delta_\mathcal{D}}(\alpha)}(\mu_{j,k})^\alpha.$$ So that, from Schwartz inequality, $$\big|\langle f,\psi_{j,k}^\lambda\rangle\big|\leq \frac{1}{\mu_{j,k}}|f|_{\Lambda_{\delta_\mathcal{D}}(\alpha)}(\mu_{j,k})^{\alpha}\Bigg( \int_{Q_{j,k}}d\mu\Bigg)\Bigg( \int_{Q_{j,k}} \psi_{j,k}^\lambda \:d\mu\Bigg)\leq |f|_{\Lambda_{\delta_\mathcal{D}}(\alpha)}(\mu_{j,k})^{\alpha+\frac{1}{2}.}$$
                
            \end{proof}
        \end{theorem}
        In order to search for reciprocals of Theorems \ref{teo61} and \ref{teo63} let us state and prove some auxiliary lemmas.
        \begin{lemma}\label{lemma64}
             Let $\mathcal{D}$ be a dyadic family in $(\mathbb{X},\mathcal{F},\mu)$ satisfying properties \ref{def21a} to \ref{def21d}. Let $\alpha>0$ be given. Assume that $f\in L^1(\mathbb{X},\mathcal{F},\mu)$ and that there exists a constant $C>0$ such that the inequality $$\big|\langle f,\psi_{j,k}^\lambda\rangle\big|\leq C(\mu_{j,k})^{\alpha+\frac{1}{2}}$$ holds for every $j \geq 0,$ every $k\in\mathcal{K}_j$ and every $\lambda=1,2,\dots,\dim\mathcal{V}_{j,k}-1.$ Then for every $Q\in\mathcal{D}$ and every $\tilde{Q}$ in the offspring of $Q$ we have the inequality $$\big|f_{\Tilde{Q}}-f_{Q}\big|\leq C\sqrt{\frac{\mu(Q)-\mu(\tilde{Q})}{\mu(Q)\mu(\tilde{Q})}}\sqrt{\dim\mathcal{V}_Q-1}\big(\mu(Q)\big)^{\alpha+\frac{1}{2}},$$ or, with index notation $$\big|f_{Q_{j+1,l}}-f_{Q_{j,k}}\big|\leq C\sqrt{\frac{\mu_{j,k}-\mu_{j+1,l}}{\mu_{j,k}\mu_{j+1,l}}}\sqrt{\dim\mathcal{V}_{j,k}-1}(\mu_{j,k})^{\alpha+\frac{1}{2}},$$ for every $l\in\mathcal{K}_{j,k}.$
             \begin{proof}
                 Let us start by obtaining a formula for $f_{\Tilde{Q}}-f_{Q}$ in terms of the Haar system $\lbrace \psi_{j,k}^\lambda\rbrace.$ Notice that from \ref{def21d}, the dimension of $\mathcal{V}_Q$ is at least two. Set $\tilde{Q}(i),$ $i=1,2,\dots,\dim\mathcal{V}_Q-1$ to denote the brothers of $\tilde{Q}.$ In other words, $Q=\tilde{Q}\cup\big(\bigcup_{i=1}^{\dim\mathcal{V}_Q-1}\tilde{Q}(i)\big),$ and the union is disjoint. Then 
                \begin{align*}
                    f_{\Tilde{Q}}-f_{Q} &= \frac{1}{\mu(\Tilde{Q})}\int_{\Tilde{Q}}  f\:d\mu-\frac{1}{\mu(Q)}\Bigg(\sum_{i=1}^{\dim\mathcal{V}_Q-1}\int_{\tilde{Q}(i)} f\:d\mu+\int_{\Tilde{Q}}f\:d\mu\Bigg) \\&=  \Bigg(\frac{1}{\mu(\Tilde{Q})}-\frac{1}{\mu(Q)}\Bigg)\int_{\Tilde{Q}}  f\:d\mu - \frac{1}{\mu(Q)}\sum_{i=1}^{\dim\mathcal{V}_Q-1}\int_{\tilde{Q}(i)} f\:d\mu\\&=  \int_{Q}\Bigg[\Bigg(\frac{1}{\mu(\Tilde{Q})}-\frac{1}{\mu(Q)}\Bigg) \chi_{\tilde{Q}} - \frac{1}{\mu(Q)}\sum_{i=1}^{\dim\mathcal{V}_Q-1}\chi_{\tilde{Q}(i)}\Bigg] f\:d\mu.
                \end{align*}
                It is clear that the function $g(x)=\Big(\frac{1}{\mu(\Tilde{Q})}-\frac{1}{\mu(Q)}\Big) \chi_{\tilde{Q}}(x) - \frac{1}{\mu(Q)}\sum_{i=1}^{\dim\mathcal{V}_Q-1}\chi_{\tilde{Q}(i)}(x)$ belongs to the space $\mathcal{V}_Q$ of functions that are constants on the offspring of $Q.$ Moreover, $$\int_\mathbb{X} g\:d\mu = \int_Q g\:d\mu = \Bigg(\frac{1}{\mu(\Tilde{Q})}-\frac{1}{\mu(Q)}\Bigg) \mu(\tilde{Q}) - \frac{1}{\mu(Q)}\big(\mu(Q)-\mu(\tilde{Q})\big)=0.$$ Hence, $g$ can be written as a linear combination of the Haar functions $\big\lbrace \psi_Q^\lambda:\lambda=1,2,\dots,\dim\mathcal{V}_Q-1\big\rbrace.$ In other words, $g=\sum_{\lambda=1}^{\dim\mathcal{V}_Q-1} C_\lambda\psi_Q^\lambda.$ Then, since $\big|\langle f,\psi_{Q}^\lambda\rangle\big|\leq C(\mu(Q))^{\alpha+\frac{1}{2}},$
                \begin{align*}
                    \big|f_{\Tilde{Q}}-f_{Q}\big| &= \Bigg|\int_\mathbb{X}\Bigg(\sum_{\lambda=1}^{\dim\mathcal{V}_Q-1} C_\lambda\psi_Q^\lambda\Bigg)f\:d\mu \Bigg| \\ &\leq \sum_{\lambda=1}^{\dim\mathcal{V}_Q-1}|C_\lambda|\big|\langle f, \psi_Q^\lambda \big\rangle\big| \\ &\leq C\Bigg(\sum_{\lambda=1}^{\dim\mathcal{V}_Q-1}|C_\lambda| \Bigg)\big(\mu(Q)\big)^{\alpha+\frac{1}{2}}.
                \end{align*}
                Let us now estimate the sum $\sum_{\lambda=1}^{\dim\mathcal{V}_Q-1}|C_\lambda|. $ From Schwartz inequality we see that $$\sum_{\lambda=1}^{\dim\mathcal{V}_Q-1}|C_\lambda| \leq \Bigg(\sum_{\lambda=1}^{\dim\mathcal{V}_Q-1}|C_\lambda|^2\Bigg)^\frac{1}{2} \sqrt{\dim\mathcal{V}_Q-1}.$$ On the other hand, since $\big\lbrace \psi_Q^\lambda:\lambda=1,2,\dots,\dim\mathcal{V}_Q-1\big\rbrace$ is orthonormal, we necessarily have that 
                \begin{align*}
                    \sum_{\lambda=1}^{\dim\mathcal{V}_Q-1}|C_\lambda|^2 &= \int_\mathbb{X} |g|^2\:d\mu \\ &= \Bigg(\frac{1}{\mu(\Tilde{Q})}-\frac{1}{\mu(Q)}\Bigg)^2 \mu(\tilde{Q}) + \frac{1}{\mu^2(Q)}\big(\mu(Q)-\mu(\tilde{Q})\big) \\ &= \frac{\mu(Q)-\mu(\tilde{Q})}{\mu(Q)\mu(\tilde{Q})}.
                \end{align*}
                Then, $$\big|f_{\Tilde{Q}}-f_{Q}\big|\leq C \sqrt{\frac{\mu(Q)-\mu(\tilde{Q})}{\mu(Q)\mu(\tilde{Q})}}\sqrt{\dim\mathcal{V}_Q-1} \big(\mu(Q)\big)^{\alpha+\frac{1}{2}}. $$
                
             \end{proof}
        \end{lemma}
        \begin{lemma}\label{lemma65}
            Let $\mathcal{D},$ $f,$ $\alpha,$ $Q$ and $\tilde{Q}$ as in Lemma \ref{lemma64}. Assume now that $\mathcal{D}$ satisfies also property \ref{def21e}. Then  $$\big|f_{\Tilde{Q}}-f_{Q}\big|\leq C \sqrt{B(B-1)}\big(\mu(Q)\big)^\alpha. $$
            \begin{proof}
                From Lemmas \ref{lemma22} and \ref{lemma64}, we have $$\big|f_{\Tilde{Q}}-f_{Q}\big|\leq C \sqrt{\frac{\mu(Q)-\mu(\tilde{Q})}{\mu(Q)\mu(\tilde{Q})}}\sqrt{B-1} \big(\mu(Q)\big)^{\alpha+\frac{1}{2}}.$$ On the other hand $$\sqrt{\frac{\mu(Q)-\mu(\tilde{Q})}{\mu(Q)\mu(\tilde{Q})}}=\sqrt{\frac{1}{\mu(\Tilde{Q})}-\frac{1}{\mu(Q)}}\leq \frac{1}{\sqrt{\mu(\tilde{Q})}}\leq \frac{\sqrt{B}}{\sqrt{\mu(Q)}},$$ and we are done.
                
            \end{proof}
        \end{lemma}
        We are now in position to prove that the power laws for the Haar coefficients of a function on $\mathbb{X}$ entail its Lipschitz regularity.
        \begin{theorem}\label{teo66}
             Let $\mathcal{D}$ be a dyadic family in $(\mathbb{X},\mathcal{F},\mu)$ satisfying properties \ref{def21a} to \ref{def21e} of Definition~\ref{def21}. Let $f$ be an integrable function in $(\mathbb{X},\mathcal{F},\mu)$ such that for some constant $C>0$ the inequalities  $$\big|\langle f,\psi_{j,k}^\lambda\rangle\big|\leq C (\mu_{j,k})^{\alpha+\frac{1}{2}}$$ hold for every $j\geq 0,$ every $k\in\mathcal{K}_j$ and every $\lambda=1,2,\dots,\dim\mathcal{V}_{j,k}-1.$ Then, with $B$ the inheritance coefficient in $\ref{def21e},$ the inequality $$\Big| f_{Q_{J,m}} - f_{Q_{J,n}}\Big|\leq \frac{2CB^\alpha\sqrt{B(B-1)}}{B^\alpha-(B-1)^\alpha} \delta_\mathcal{D}^\alpha(Q_{J,m},Q_{J,n})$$ holds for every positive integer $J$ and every $m,n\in\mathcal{K}_J.$ In other words, $|f|_{\Lambda_{\delta_\mathcal{D}(\alpha)}}\leq \frac{2CB^\alpha\sqrt{B(B-1)}}{B^\alpha-(B-1)^\alpha}.$
             \begin{proof}
                 Fix $J,$ $m$ and $n.$ Of course, we may assume that $m\neq n.$ Let $Q_{j,k}\in\mathcal{D}$ be the smallest common ancestor of both $Q_{J,m}$ and $Q_{J,n}.$ Observe that, since $\mathbb{X}$ itself is a common ancestor of $Q_{J,m}$ and $Q_{J,n},$ such a $Q_{j,k}$ is well-defined. Moreover, $\mu_{j,k}=\mu(Q_{j,k})=\delta_\mathcal{D}(Q_{J,m},Q_{J,n}).$ It is clear that the scale index $j$ of $Q_{j,k}$ is smaller than $J.$ Let $l$ be the positive integer such that $j+l=J.$ Then, we may follow the genealogies of $Q_{J,m}$ and $Q_{J,n}$ up to their common ancestor $Q_{j,k}$ in the following way $$ Q_{J,m}\subset Q_{J-1,k(J-1,J,m)} \subset Q_{J-2,k(J-2,J,m)} \subset \dots \subset Q_{j+1,k(j+1,J,m)} \subset Q_{j,k},$$ and $$ Q_{J,n}\subset Q_{J-1,k(J-1,J,n)} \subset Q_{J-2,k(J-2,J,n)} \subset \dots \subset Q_{j+1,k(j+1,J,n)} \subset Q_{j,k},$$ where the notation $Q_{J-i,k(J-i,J,m)}$ stands for $i^{th}$ ancestor of $Q_{J,m}.$ With these two family chains we may write
                \begin{align*}
                    f_{Q_{J,m}}-f_{Q_{J,n}} &= \big(f_{Q_{J,m}}- f_{ Q_{j,k}}\big) + \big(f_{Q_{j,k}} - f_{ Q_{J,n}}\big)\\ &= \sum_{i=0}^{l-1}\big(f_{Q_{J-i,k(J-i,J,m)}} - f_{Q_{J-(i+1),k(J-(i+1),J,m)}} \big) \\ &\:\:\: - \sum_{i=0}^{l-1}\big(f_{Q_{J-i,k(J-i,J,n)}} - f_{Q_{J-(i+1),k(J-(i+1),J,n)}} \big).
                \end{align*}
                Now, taking absolute values and using Lemma \ref{lemma65} in every term of the sums we obtain 
                 \begin{align*}
                    \big|f_{Q_{J,m}}-f_{Q_{J,n}} \big|&\leq \sum_{i=0}^{l-1}\big|f_{Q_{J-i,k(J-i,J,m)}} - f_{Q_{J-(i+1),k(J-(i+1),J,m)}} \big| \\ &\:\:\: - \sum_{i=0}^{l-1}\big|f_{Q_{J-i,k(J-i,J,n)}} - f_{Q_{J-(i+1),k(J-(i+1),J,n)}} \big| \\ &\leq C\sqrt{B(B-1)}\sum_{i=0}^{l-1} \big(\mu(Q_{J-(i+1),k(J-(i+1),J,m)})\big)^\alpha \\ &\:\:\: + C\sqrt{B(B-1)}\sum_{i=0}^{l-1}\big(\mu(Q_{J-(i+1),k(J-(i+1),J,n)})\big)^\alpha.
                \end{align*}
                Since both subcubes sequences are contained in and increasing to $Q_{j,k},$ we apply \ref{lem223} in Lemma \ref{lemma22} to get 
                \begin{align*}
                    \big|f_{Q_{J,m}}-f_{Q_{J,n}}\big| &\leq 2C\sqrt{B(B-1)}\Bigg( \sum_{i=0}^{l-1}\Big({1-\frac{1}{B}}\Big)^{(l-i-1)\alpha}\Bigg)(\mu_{j,k})^\alpha \\ &\leq \frac{2CB^\alpha\sqrt{B(B-1)}}{B^\alpha-(B-1)^\alpha}\delta_\mathcal{D}^\alpha(Q_{J,m},Q_{J,n}).
                \end{align*}

             \end{proof}
        \end{theorem}
        
        We have already observed that our general dyadic families with inheritance coefficient $B$ need not be differentiation bases for $L^1(\mathbb{X},\mathcal{F},\mu)$ functions. Nevertheless, the results above show that a power law in the scale size of the Haar coefficients, are still giving some information on the smoothness of the means of a function over the dyadic sets in the given family. To be a differentiation basis for a dyadic family as defined above means that for every $f\in L^1(\mathbb{X},\mathcal{F},\mu),$ $$ \lim_{j\to\infty} \frac{1}{\mu(Q_{j,k(j,x)})}\int_{Q_{j,k(j,x)}}f(y)\:d\mu(y)=f(x)$$ for $\mu$ almost every $x\in\mathbb{X}.$ Since, as it easy to prove, the dyadic maximal function $$\mathcal{M}_{\mathcal{D}} f(x)= \sup\Bigg\lbrace\frac{1}{\mu(Q)}\int_Q |f|\:d\mu \Bigg\rbrace$$ is of weak type $(1,1),$ what is generally lacking in our setting is the existence of some dense subspace of $L^1(\mathbb{X},\mathcal{F},\mu)$ for which the limit above is satisfied. Of course, we recover this properties in general settings such as spaces of homogeneous type with metric guided constructions of dyadic families, such as those in \cite{Chr90}. The next result is a consequence of Theorem \ref{teo66} under the hypothesis of differentiation and can be seen as a reciprocal of Theorem \ref{teo61} above.

        \begin{corollary}\label{coro67}
            Let $\mathcal{D}$ be a dyadic family in $(\mathbb{X},\mathcal{F},\mu)$ satisfying properties \ref{def21a} to \ref{def21e}. Assume also that $\mathcal{D}$ is a differentiation basis. Let $f$ be an integrable function in $(\mathbb{X},\mathcal{F},\mu)$ such that for some constant $C>0$ the inequalities $$ \big|\langle f,\psi_{j,k}^\lambda\rangle\big|\leq C (\mu_{j,k})^{\alpha+\frac{1}{2}}$$ holds for every $j\geq 0$, every $k\in\mathcal{K}_j$ and every $\lambda=1,2,\dots,\dim\mathcal{V}_{j,k}-1.$ then, with $B$ the inheritance coefficient in \ref{def21e}, we have the inequality $$\big| f(x)-f(y)\big|\leq \frac{2CB^\alpha\sqrt{B(B-1)}}{B^\alpha-(B-1)^\alpha} \delta_\mathcal{D}^\alpha(x,y),$$ for almost every $x$ and $y$ in $\mathbb{X}.$ Moreover $|f|_{Lip_\mathcal{D}(\alpha)}\leq \frac{2CB^\alpha\sqrt{B(B-1)}}{B^\alpha-(B-1)^\alpha}.$
            \begin{proof}
                Given $x$ and $y$ in $\mathbb{X},$ let us consider the two sequences of cubes $\lbrace Q_{j,k(j,x)}:j \geq 0 \rbrace$ and $\lbrace Q_{j,k(j,y)}:j \geq 0 \rbrace,$ and let us apply the result of Theorem \ref{teo66} to these sequences, hence $$\Big| f_{Q_{j,k(j,x)}} - f_{Q_{j,k(j,y)}}\Big|\leq \frac{2CB^\alpha\sqrt{B(B-1)}}{B^\alpha-(B-1)^\alpha} \delta_\mathcal{D}^\alpha(Q_{j,k(j,x)},Q_{j,k(j,y)}).$$ Because of the differentiation property of $\mathcal{D}$ for $L^1(\mathbb{X},\mathcal{F},\mu)$ functions, the left-hand side of the above inequality tends to $\big| f(x)-f(y)\big|$ when $j$ grows to infinity for almost every $x$ and almost every $y,$ both in $\mathbb{X}.$ On the other hand, let $Q_{j_0,k_0}$ be the smallest dyadic cube containing both $x$ and $y.$ Then $\delta_{\mathcal{D}}(x,y)=\mu_{j_0,k_0}.$ Since for $j>j_0+1$ we have that $Q_{j,k(j,x)}$ and $Q_{j,k(j,y)}$ are disjoint subcubes of $Q_{j_0,k_0},$ we necessarily have that $$\delta_{\mathcal{D}}(Q_{j,k(j,x)},Q_{j,k(j,y)})\leq \mu_{j_0,k_0}=\delta_{\mathcal{D}}(x,y),$$ and we are done.
                
            \end{proof}
        \end{corollary}

    \section{Examples} \label{sec:Section7}
        The basic results in the above section relate the dyadic regularity of a function on $\mathbb{X}$ with a power law for the Haar coefficients in terms of the measures $\mu_{j,k}$ of the supports of the $\psi_{j,k}^\lambda.$ One may take averages on $k\in\mathcal{K}_j$ and $\lambda\in\lbrace1,\dots, \dim \mathcal{V}_{j,k}-1\rbrace$ for fixed $j$ in the inequalities for the wavelet coefficients provided by Theorem \ref{teo66} and Corollary \ref{coro67};  $ \big|\langle f,\psi_{j,k}^\lambda\rangle\big|\leq C (\mu_{j,k})^{\alpha+\frac{1}{2}},$ $j\geq 0,$ $k\in\mathcal{K}_j$ and $\lambda=1,2,\dots,\dim\mathcal{V}_{j,k}-1.$ In this way we obtain the average wavelet coefficients $$AWC(f)_{(j)}\doteq\frac{1}{K_j}\sum_{k\in\mathcal{K}_j}\frac{1}{\dim\mathcal{V}_{j,k}-1}\sum_{\lambda=1}^{\dim\mathcal{V}_{j,k}-1}\big|\langle f,\psi_{j,k}^\lambda\rangle\big|\leq C \frac{1}{K_j} \sum_{k\in\mathcal{K}_j}(\mu_{j,k})^{\alpha+\frac{1}{2}}.$$
        When $\mu_{j,k}=\mu_j$ only depends on the scale $j$ but not on the position $k$ of the support of $\psi_{j,k}^\lambda,$ we have that $AWC(f)_{(j)}\leq C (\mu_{j})^{\alpha+\frac{1}{2}}$ for every $j\geq 0.$ This functional inequality, involving two parameters $\alpha$ and $C,$ the Lipschitz exponent and the corresponding bound according to Definition~\ref{def44},  can be written in simple form by taking logarithms, $$LAWC(f)_{(j)}\doteq \log AWC(f)_{(j)} \leq \Big( \alpha+\frac{1}{2}\Big)\log \mu_j + \log C.$$ The important information in this inequality is that for empirical images $AWC(f)_{(j)}$ and $LAWC(f)_{(j)}$ can explicitly be computed. Sometimes, in particular for uniform partitions, $LAWC$ behaves like a linear function of $j$ with slope $\beta<0,$ and this fact gives a plausible value of the regularity exponent $\alpha=-\beta-\frac{1}{2}.$ 
    
\begin{figure}[H]
            \centering
            \begin{subfigure}{0.19\textwidth}
                \centering
                \includegraphics[width=\textwidth]{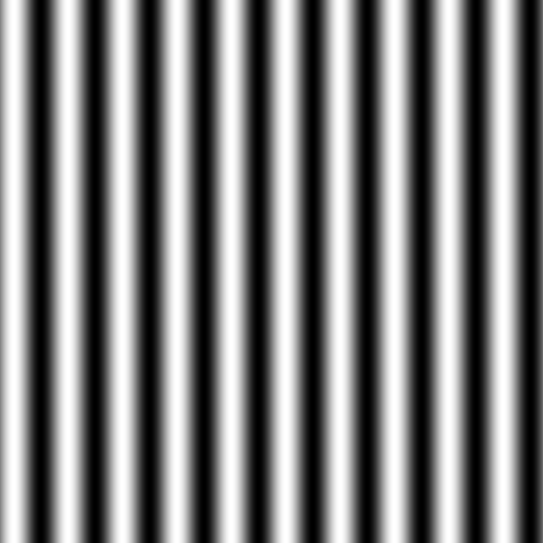}
                \caption*{$F_1$}
                \label{fig:figura1}
            \end{subfigure}
            \hspace{0.5cm}
            \begin{subfigure}{0.19\textwidth}
                \centering
                \includegraphics[width=\textwidth]{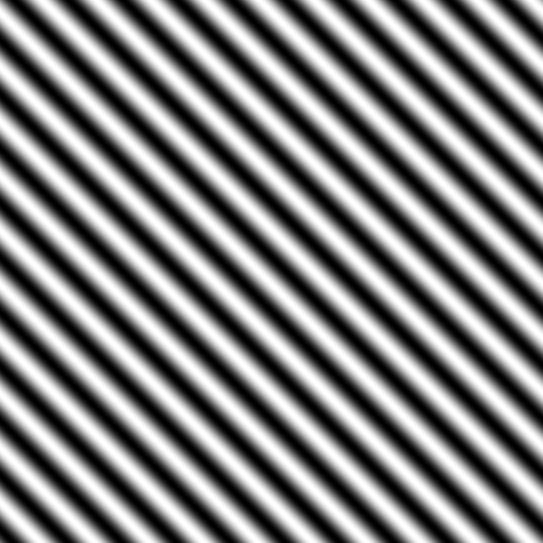}
                \caption*{$F_2$}
                \label{fig:figura2}
            \end{subfigure}
            
            \caption{The level sets of $F_1$ and $F_2$. Black corresponds to the value $0$ and white to $255.$}
            \label{fig:figure1}
        \end{figure}

\begin{figure}[H]
            \centering
            \begin{subfigure}{0.19\textwidth}
                \centering
                \includegraphics[width=\textwidth]{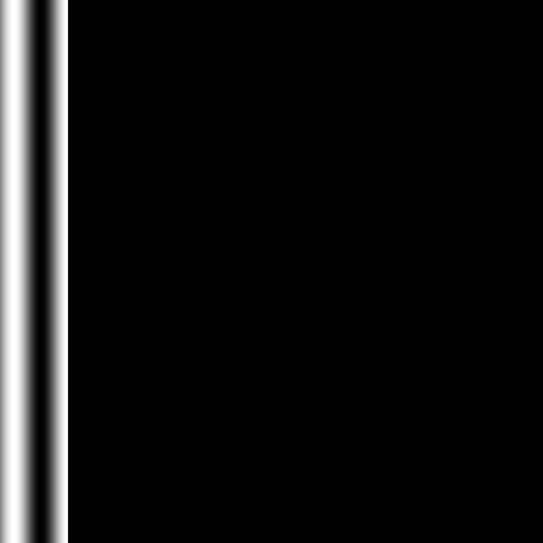}
                \caption*{$G_1$}
                \label{fig:figura7}
            \end{subfigure}
            \hfill
            \begin{subfigure}{0.19\textwidth}
                \centering
                \includegraphics[width=\textwidth]{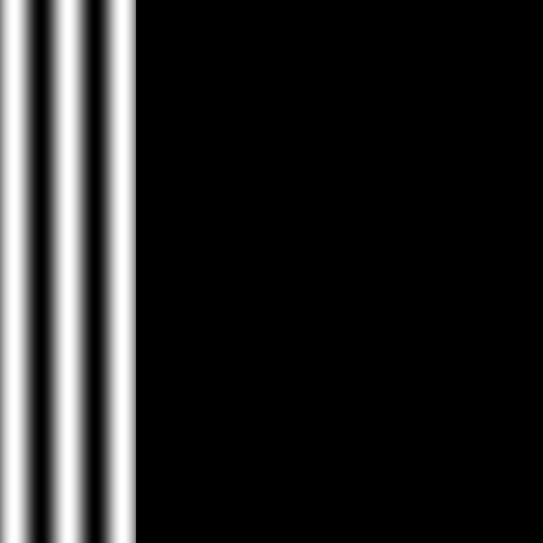}
                \caption*{$G_2$}
                \label{fig:figura6}
            \end{subfigure}
            \hfill
            \begin{subfigure}{0.19\textwidth}
                \centering
                \includegraphics[width=\textwidth]{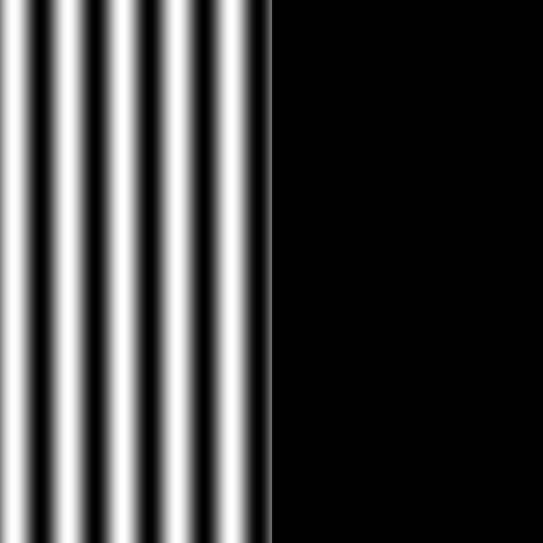}
                \caption*{$G_3$}
                \label{fig:figura8}
            \end{subfigure}
            \hfill
             \begin{subfigure}{0.19\textwidth}
                \centering
                \includegraphics[width=\textwidth]{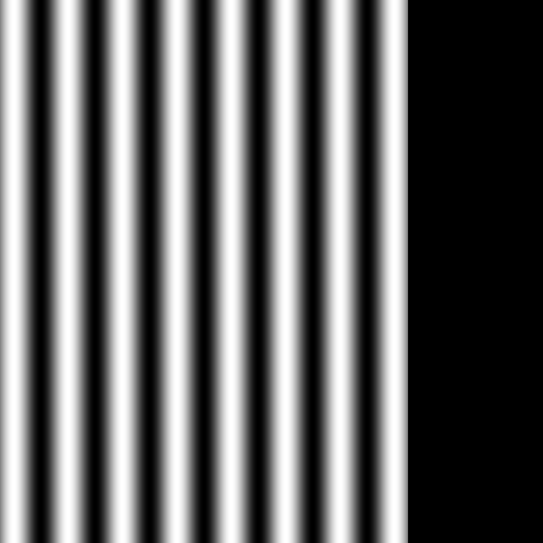}
                \caption*{$G_4$}
                \label{fig:figura9}
            \end{subfigure}
            \hfill
             \begin{subfigure}{0.19\textwidth}
                \centering
                \includegraphics[width=\textwidth]{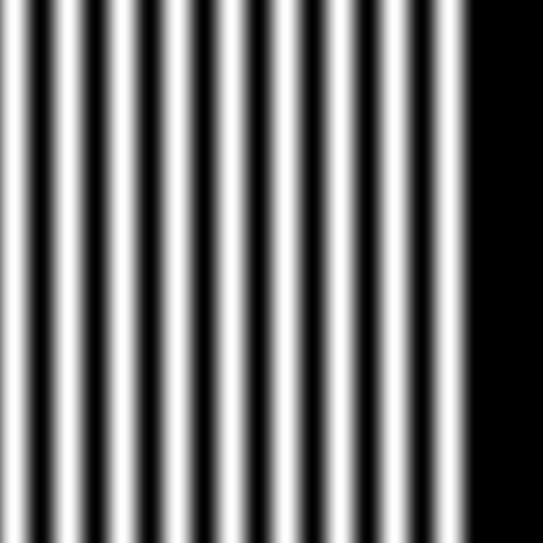}
                \caption*{$G_5$}
                \label{fig:figura10}
            \end{subfigure}
            \caption{The five functions $G_i,$ $ i=1,\dots,5.$}
            \label{fig:figure2}
        \end{figure}

        In this section, we aim to illustrate, using some simple two-dimensional images, the potential use of diverse dyadic systems to detect particular textures in images. By linear regression, we shall estimate, using Python, the slope and the intercept of the line, as a function of the scale $j,$ on the upper estimate of $LAWC(f)_{(j)}.$ The slope $-(\beta+\frac{1}{2})$ gives an estimate of the worst regularity in the image, and the intercept $\log C$ provides an estimate for the amount of that regularity. Let us precise the images and the dyadic systems that we shall consider. The first set of images contains the two following trigonometric bivariate functions on the square $[0,1]^2=[0,1]\times [0,1]$ with the usual area measure, $F_1(x,y)=\frac{255}{2} \big(1 + \sin(20\pi x) \big)$ and
        $F_2(x,y)=\frac{255}{2}  \big(1 + \sin(20\pi (x+y)) \big).$ The second set of images on $[0,1]^2$ is given by $G_i(x,y) = F_1(x,y)\chi_{[0,b_i]} (x,y),$ $i=1,\dots,5;$ with $b_1= \frac{1}{8},$ $b_2= \frac{1}{4},$ $b_3= \frac{1}{2},$ $b_4= \frac{3}{4}$ and $b_5= \frac{7}{8}.$ The first set $\lbrace F_1,F_2\rbrace$ is depicted in Figure~\ref{fig:figure1} and the second in Figure~\ref{fig:figure2}.

        The three families of wavelets that we shall use are associated with the dyadic families that we introduced in Section~\ref{sec:Section2}. Namely, squares, parabolic rectangles, and non-metric bands. Precisely, we shall get the Haar coefficients of $F_i,$ $i=1,2$ and $G_i,$ $i=1,\dots, 5$ with respect to three wavelets, one in each of the dyadic cases. The first is the classical Haar function of the plane, given by       $$\psi_{j,\mathbf{k}}^{\mathrm{I}}(x,y)=2^j h_{0,0}(2^jx-k_1)\chi_{[0,1]}(2^jy-k_2),$$
        with $\mathbf{k}=(k_1,k_2),$ $0\leq k_i\leq 2^j-1,$ $i=1,2$ and $h_{0,0}(t)=\chi_{[0,\frac{1}{2})}(t)-\chi_{[\frac{1}{2},1)}(t)$ the one-dimensional Haar function. The second is one of the seven $(2^3-1)$ Haar parabolic functions given by 
        $$\psi_{j,\mathbf{k}}^{\mathrm{II}}(x,y)=8^\frac{j}{2} h_{0,0}(4^jx-k_1)\chi_{[0,1]}(2^jy-k_2),$$ with $\mathbf{k}=(k_1,k_2),$ $0\leq k_1\leq 4^j-1$ and $0\leq k_2\leq 2^j-1.$ The third is the associated non-differentiating dyadic family of vertical bands,
        $$\psi_{j,k}^{\mathrm{III}}(x,y)=2^\frac{j}{2} h_{0,0}(2^jx-k),$$ with $0 \leq k\leq 2^j-1.$ Each one of these three wavelets are applied to obtain the corresponding Haar coefficients of each one of the seven images. In this way, we obtain for $j=0,\dots,12,$ the values of $LAWC(F_i, \mathrm{I})_{(j)},$ $LAWC(F_i, \mathrm{II})_{(j)}$ and $LAWC(F_i, \mathrm{III})_{(j)}$ for $i=1,2$ and $LAWC(G_i, \mathrm{I})_{(j)},$ $LAWC(G_i, \mathrm{II})_{(j)}$ and $LAWC(G_i, \mathrm{III})_{(j)}$ for $i=1,\dots,5.$ The estimates of the Lipschitz exponents for $F_1$ and $F_2$ are summarized in Table~\ref{tab:table1}.
        
        \renewcommand{\arraystretch}{1.} % increases the high of rows
        \begin{table}[!h]
            \centering
            \begin{tabularx}{\textwidth}{|>{\centering\arraybackslash}m{1.cm}|>{\centering\arraybackslash}X|>{\centering\arraybackslash}X|}
                \rowcolor{lightgray}
                \hline 
                & $F_1$ & $F_2$ \\ \hline 
                \cellcolor{lightgray} $\mathrm{I}$ & $\alpha_{1,\mathrm{I}}= 0.4936$ & $\alpha_{2,\mathrm{I}}= 0.4782$ \\ \hline
                \cellcolor{lightgray} $\mathrm{II}$ & $\alpha_{1,\mathrm{II}}= 0.6234$ & $\alpha_{2,\mathrm{II}}= 0.3002$ \\ \hline
                \cellcolor{lightgray} $\mathrm{III}$ & $\alpha_{1,\mathrm{III}}= 0.4873$ & $\alpha_{2,\mathrm{III}}= 0.3412$ \\ \hline
            \end{tabularx}
            \caption{The three exponents of regularity for $F_1$ and $F_2$.}
            \label{tab:table1}
        \end{table}
        
        The estimates for the intercepts of the regression lines for $G_i,$ $i=1,\dots,5$ are summarized in Table~\ref{tab:table2}, where aside from the increment of $\log C$ for $i$ increasing in the set $\lbrace1,2,3,4,5\rbrace,$ shows the empirical stability of the regularity exponent $\alpha.$ To make the comparison possible, we use different bases for the logarithm function according to the scales of the dyadic objects.

        \begin{table}[!h]
            \centering
            \begin{tabularx}{\textwidth}{|>{\centering\arraybackslash}m{1cm}
                            |>{\centering\arraybackslash}X
                            |>{\centering\arraybackslash}X
                            |>{\centering\arraybackslash}X
                            |>{\centering\arraybackslash}X
                            |>{\centering\arraybackslash}X|}
                \rowcolor{lightgray}
                \hline
                     & {$G_1$} & {$G_2$} & {$G_3$} & {$G_4$} & {$G_5$} \\ \hline \cellcolor{lightgray}
                     $\mathrm{I}$ &  $\alpha_{1,\mathrm{I}} =0.4491$ 
                    
                     $\log_4 C_{1,\mathrm{I}}=4.8237$
                     &  $\alpha_{2,\mathrm{I}}=0.4821$
        
                     $\log_4 C_{2,\mathrm{I}}=5.5285$
                     & $\alpha_{3,\mathrm{I}}=0.4892$
        
                     $\log_4 C_{3,\mathrm{I}}=6.0666$
                     & $\alpha_{4,\mathrm{I}}=0.4896$ 
                     
                     $\log_4 C_{4,\mathrm{I}}=6.3640$& 
                     $\alpha_{5,\mathrm{I}}=0.4933$ 
                     
                     $\log_4 C_{5,\mathrm{I}}=6.4899$\\ 
                     \hline 
                     \cellcolor{lightgray}
                     $\mathrm{II}$ & $\alpha_{1,\mathrm{II}}=0.5640$ 
                    
                     $\log_8 C_{1,\mathrm{II}}=3.6739$
                     &  $\alpha_{2,\mathrm{II}}=0.6087$
        
                     $\log_8 C_{2,\mathrm{II}}=4.1613$
                     & $\alpha_{3,\mathrm{II}}=0.6174$
        
                     $\log_8 C_{3,\mathrm{II}}=4.5212$
                     & $\alpha_{4,\mathrm{II}}=0.6176$ 
                     
                     $\log_8 C_{4,\mathrm{II}}=4.7176$& 
                     $\alpha_{5,\mathrm{II}}=0.6207$ 
                     
                     $\log_8 C_{5,\mathrm{II}}=4.7941$ \\ \hline \cellcolor{lightgray}
                    $\mathrm{III}$ &  $\alpha_{1,\mathrm{III}}=0.3981$ 
                    
                     $\log_2 C_{1,\mathrm{III}}=8.6474$
                     &  $\alpha_{2,\mathrm{III}}=0.4642$
        
                     $\log_2 C_{2,\mathrm{III}}=10.0569$
                     & $\alpha_{3,\mathrm{III}}=0.4784$
        
                     $\log_2 C_{3,\mathrm{III}}=11.1332$
                     & $\alpha_{4,\mathrm{III}}=0.4793$ 
                     
                     $\log_2 C_{4,\mathrm{III}}=11.7281$&
                     $\alpha_{5,\mathrm{III}}=0.4865$ 
                     
                     $\log_2 C_{5,\mathrm{III}}=11.9798$ \\ \hline
                    
            \end{tabularx}
            \caption{Persistence of $\alpha$ and the increment of $\log C.$}
            \label{tab:table2}
        \end{table}
        
        Notice that, from Table~\ref{tab:table1}, the three types of wavelets are detecting more regularity in $F_1$ than in $F_2,$ but the parabolic wavelets have a better separation performance between the vertical and oblique textures in these images.

    %\input{referencias.bbl}  
    %\bibliographystyle{plain}
    %\bibliography{referencias}
%
%

%%%%%%%%%%%%%  References %%%%%%%%%%%%%%%%%
%\bibliographystyle{amsalpha}
%\bibliography{ref}

%\nocite{*}

%%%%%%%%%%%%%%%%%%%%% Declarations %%%%%%%%%%%%%%%%%%%%%%%%%%%%%%%%%%%%%%

%\section*{Declarations}
%
%\subsection*{Acknowledgements}
%This work was supported by Consejo Nacional de Investigaciones Cient\'ificas y T\'ecnicas - CONICET and Universidad Nacional del Litoral - UNL in Argentina.
%
%\subsection*{Funding}
%Consejo Nacional de Investigaciones Cient\'ificas y T\'ecnicas, grant PIP-2021-2023-11220200101940CO.
%
%\subsection*{Conflicts of interest/Competing interests}
%The authors have no conflicts of interest to declare that are relevant to the content of this article.
%
%\subsection*{Availability of data and material}
%Not applicable.

%%%%%%%%%%%%% Address Authors %%%%%%%%%%%%%%%

\bigskip

\medskip

\noindent{\footnotesize
\noindent\textit{Affiliation.\,}
\textsc{Instituto de Matem\'{a}tica Aplicada del Litoral ``Dra. Eleonor Harboure'', UNL, CONICET.}
	
\smallskip
\noindent\textit{Address.\,}\textmd{IMAL, Streets F.~Leloir and A.P.~Calder\'on, CCT CONICET Santa Fe, Predio ``Alberto Cassano'', Colectora Ruta Nac.~168 km~0, Paraje El Pozo, S3007ABA Santa Fe, Argentina.}
	
\smallskip

}

\bigskip

\end{document}